\documentclass{gtmon_a}
\pdfoutput=1

\usepackage{amscd}

\proceedingstitle{Groups, homotopy and configuration spaces (Tokyo
  2005)}
\conferencestart{5 July 2005}
\conferenceend{11 July 2005}
\conferencename{Groups, homotopy and configuration spaces, 
                in honour of Fred Cohen's 60th birthday}
\conferencelocation{University of Tokyo, Japan}

\editor{Norio Iwase}
\givenname{Norio}
\surname{Iwase}

\editor{Toshitake Kohno}
\givenname{Toshitake}
\surname{Kohno}

\editor{Ran Levi}
\givenname{Ran}
\surname{Levi}

\editor{Dai Tamaki}
\givenname{Dai}
\surname{Tamaki}

\editor{Jie Wu}
\givenname{Jie}
\surname{Wu}

\title [Basis-conjugating automorphisms of a free group]
{Basis-conjugating automorphisms of a free group\\
and associated Lie algebras}

\author{F\,R Cohen}
\givenname{F\,R}
\surname{Cohen}
\email{cohf@math.rochester.edu}

\author{J Pakianathan}
\givenname{J}
\surname{Pakianathan}
\email{jonpak@math.rochester.edu}

\author{V\,V Vershinin}
\givenname{V\,V}
\surname{Vershinin}
\address{{\rm FRC and JP:}\ \ Department of Mathematics\\
University of Rochester\\\newline
Rochester, NY 14627\\USA\vspace{3pt}\\\newline
{\rm VVV:}\ \ D\'epartement des Sciences Math\'ematiques\\
Universit\'e Montpellier II\\\newline
Place Eug\'ene Bataillon\\34095
Montpellier cedex 5\\France\\\newline
and\\\newline
Sobolev Institute of Mathematics\\Novosibirsk, 630090\\Russia\vspace{3pt}\\\newline
{\rm JW:}\ \ Department of Mathematics\\
National University of Singapore\\
Singapore 117543\\
Republic of Singapore}
\email{vershini@math.univ-montp2.fr}
\email{versh@math.nsc.ru}

\author{J Wu}
\givenname{Jie}
\surname{Wu}
\email {matwuj@nus.edu.sg}

\volumenumber{13}
\issuenumber{}
\publicationyear{2008}
\papernumber{06}
\startpage{147}
\endpage{168}

\doi{}
\MR{}
\Zbl{}

\arxivreference{math.AT/0610946} 


\received{31 May 2006}
\revised{6 March 2007}
\accepted{13 March 2007}
\published{22 February 2008}
\publishedonline{22 February 2008}
\proposed{}
\seconded{}
\corresponding{}
\version{}

\makeatletter
\def\cnewtheorem#1[#2]#3{\newtheorem{#1}{#3}[section]
\expandafter\let\csname c@#1\endcsname\c@thm}

\makeop{Aut}
\let\xysavmatrix\xymatrix
\def\xymatrix{\disablesubscriptcorrection\xysavmatrix}

\newtheorem{thm}{Theorem}[section]
\cnewtheorem{cor}[thm]{Corollary}
\cnewtheorem{lem}[thm]{Lemma}
\cnewtheorem{prop}[thm]{Proposition}

\theoremstyle{definition}
\cnewtheorem{defin}[thm]{Definition}
\theoremstyle{definition}
\cnewtheorem{exm}[thm]{Example}

\theoremstyle{remark}

\makeatother  %  move after \newtheorem block

\begin{document}

\begin{htmlabstract}
Let F<sub>n</sub> = &lang;x<sub>1</sub>, &hellip;, x<sub>n</sub>&rang;
denote the free group with generators {x<sub>1</sub>, &hellip;,
x<sub>n</sub>}. Nielsen and Magnus described generators for the kernel of
the canonical epimorphism from the automorphism group of F<sub>n</sub>
to the general linear group over the integers. In particular among them
are the automorphisms &chi;<sub>k,i</sub> which conjugate the generator
x<sub>k</sub> by the generator x<sub>i</sub> leaving the x<sub>j</sub>
fixed for j not k. A computation of the cohomology ring as well as
the Lie algebra obtained from the descending central series of the
group generated by &chi;<sub>k,i</sub> for i < k is given here. Partial
results are obtained for the group generated by all &chi;<sub>k,i</sub>.
\end{htmlabstract}

\begin{abstract}
Let $F_n = \langle x_1, \dots, x_n\rangle$ denote the free group with
generators $\{x_1, \dots, x_n\}$. Nielsen and Magnus described
generators for the kernel of the canonical epimorphism from the
automorphism group of $F_n$ to the general linear group over the
integers. In particular among them are the automorphisms $\chi_{k,i}$
which conjugate the generator $x_k$ by the generator $x_i$ leaving the
$x_j$ fixed for $j \neq k$. A computation of the cohomology ring as
well as the Lie algebra obtained from the descending central series of
the group generated by $\chi_{k,i}$ for $i<k$ is given here. Partial
results are obtained for the group generated by all $\chi_{k,i}$.
\end{abstract}

\begin{asciiabstract} 
Let F_n = <x_1,...,x_n> denote the free group with generators
{x_1,...,x_n}.  Nielsen and Magnus described generators for the
kernel of the canonical epimorphism from the automorphism group of F_n
to the general linear group over the integers.  In particular among
them are the automorphisms chi_{k,i} which conjugate the generator x_k
by the generator x_i leaving the x_j fixed for j not k.  A computation
of the cohomology ring as well as the Lie algebra obtained from the
descending central series of the group generated by chi_{k,i} for i<k
is given here. Partial results are obtained for the group generated by
all chi_{k,i}.
\end{asciiabstract}

\maketitle

\section{Introduction}
Let $\pi$ be a discrete group with $\Aut(\pi)$ the automorphism group
of $\pi$. Consider the free group $F_n$ generated by $n$ letters
$\{x_1, x_2, \cdots, x_n\}$. The kernel of the natural map
$$\Aut(F_n) \to GL(n, \mathbb Z)$$ is denoted $IA_n$. Nielsen, and
Magnus gave automorphisms which generate $IA_n$ as a group
\cite{N,magnus,mks}. These automorphisms are named as follows:

\begin{itemize}
        \item $\chi_{k,i}$ for $ i \neq k$ with $1 \leq i,k \leq n$, and
        \item $\theta (k;[s,t])$ for $k,s,t$ distinct integers with $1 \leq k,s,t \leq n$ and $s<t$.
\end{itemize}

The definition of the map $\chi_{k,i}$ is given by the formula
\[
\chi_{k,i}(x_j)=
\begin{cases}
x_j & \text{if $k \neq j$,}\\
(x_i^{-1})(x_k)(x_i) & \text{if $k = j$.}
\end{cases}
\] Thus the map $\chi_{k,i}^{-1}$ satisfies the formula
\[
\chi_{k,i}^{-1}(x_j)=
\begin{cases}
x_j & \text{if $k \neq j$,}\\
(x_i)(x_k)({x_i}^{-1}) & \text{if $k = j$.}
\end{cases}
\]
The map $\theta(k;[s,t])$ is defined by the formula
\[
\theta(k;[s,t])(x_j)=
\begin{cases}
x_j & \text{if $k \neq j$,}\\
(x_k)\cdot([x_s,x_t]) & \text{if $k = j$.}
\end{cases}
\] for which the commutator is given by 
$[a,b] = a^{-1}\cdot b^{-1}\cdot a \cdot b$. Thus
the map $\theta(k;[s,t])^{-1}$ satisfies the formula
\[
\theta(k;[s,t])^{-1}(x_j)=
\begin{cases}
x_j & \text{if $k \neq j$,}\\
(x_k)\cdot([x_s,x_t]^{-1}) & \text{if $k = j$.}
\end{cases}
\]
Consider the subgroup of $IA_n$ generated by the $\chi_{k,i}$, the
group of basis conjugating automorphisms of a free group. This
subgroup has topological interpretations. First of all it is the pure
group of motions of $n$ unlinked circles in $S^3$ (Goldsmith \cite{dg}
and Jensen, McCammond and Meier \cite{jmm}) and because of this it is
known as the ``group of loops". On the other hand it is also the pure
braid--permutation group. This is explained at the end of this
section. This group is denoted $P\varSigma_n$ in \cite{jmm}. McCool
gave a presentation for it \cite{mc}. This presentation is listed in
\fullref{thm:McCool relations} below. The subgroup of $P\varSigma_n$
generated by the $\chi_{k,i}$ for $ i < k $ is denoted
$P\varSigma_n^+$ here and is called the ``upper triangular McCool
group" in \cite{ccp}.

The purpose of this article is to determine the natural Lie algebra
structure obtained from the descending central series for
$P\varSigma_n^+$ together with related information for $P\varSigma_n$
as well as the structure of the cohomology ring of
$P\varSigma_n^+$. One motivation for the work here is that the groups
$P\varSigma_n$ and $P\varSigma_n^+$ are natural as well as accessible
cases arising as analogues of work in D~Johnson \cite{j}, S~Morita
\cite{m}, D~Hain \cite{h}, B~Farb \cite{f}, N~Kawazumi \cite{ka},
T~Kohno \cite{k}, C~Jensen, J~McCammond, and J~Meier \cite{jmm},
T~Sakasai \cite{s}, T~Satoh \cite{sa}, A~Pettet \cite{pe} and Y~Ihara
\cite{y}.  In those works the Johnson filtration is used frequently
rather than the descending central series. The techniques here for
addressing these Lie algebras are due to T~Kohno \cite{k} and M~Falk
and R~Randell \cite{fr}.

The cohomology of $P\varSigma_n$ was computed by C~Jensen,
J~McCammond, and J~Meier \cite{jmm}.  N~Kawazumi \cite{ka},
T~Sakasai \cite{s}, T~Satoh \cite{sa} and A~Pettet \cite{pe} have
given related cohomological information for $IA_n$. The integral
cohomology of the natural direct limit of the groups $\Aut(F_n)$ is
given in work of S~Galatius \cite{g}.

The main results here arise from McCool's presentation which is
stated next.

\begin{thm} \label{thm:McCool relations}
A presentation of $P\varSigma_n$ is given by generators $\chi_{k,j}$
together with the following relations.
\begin{enumerate}
\item $\chi_{i,j}\cdot \chi_{k,j}\cdot \chi_{i,k} = \chi_{i,k} \cdot
\chi_{i,j} \cdot \chi_{k,j}$ for $i,j,k$ distinct.

\item $[\chi_{k,j}, \chi_{s,t}] = 1$ if $\{j,k\} \cap \{s,t \}=
\emptyset$.

\item $[ \chi_{i,j},\chi_{k,j}] = 1$ for
$i,j,k$ distinct.

\item $[ \chi_{i,j} \cdot \chi_{k,j}, \chi_{i,k}] = 1$ for
$i,j,k$ distinct (redundantly).

\end{enumerate}
\end{thm}

In what follows below, $gr^*(\pi)$ denotes the associated graded Lie
algebra obtained from the descending central series of a discrete
group $\pi$. Work of T~Kohno \cite{k}, as well as M~Falk, and
R~Randell \cite{fr} provide an important description of these Lie
algebras for many groups $\pi$, one of which is the pure braid group
on $n$ strands $P_n$. The Lie algebra $gr^*(P_n)$, basic in Kohno's
work, gave an important ingredient in his analysis of Vassiliev
invariants of pure braids in terms of iterated integrals \cite{k}. A
presentation for this Lie algebra is given by the quotient of the
free Lie algebra $L[B_{i,j}| 1 \leq i < j \leq k]$ generated by
elements $B_{i,j}$ with $1 \leq i < j \leq k$ modulo the
``infinitesimal braid relations" or ``horizontal $4T$ relations"
given by the following three relations:

\begin{enumerate}
  \item If $\{i,j\} \cap \{s,t\} = \emptyset$, then
  $[B_{i,j}, B_{s,t}] = 0$.
  \item If $ i<j<k$, then $[B_{i,j}, B_{i,k} + B_{j,k}] = 0$.
  \item If $ i<j<k$, then $[B_{i,k}, B_{i,j} + B_{j,k}] = 0$.
\end{enumerate}

The results below use the methods of Kohno, and Falk--Randell to
obtain information about $P\varSigma_n$, and $P\varSigma_n^+$ as
well. One feature is that the Lie algebras given by
$gr^*(P\varSigma_n)$ and $gr^*(P\varSigma_n^+)$ satisfy two of the
``horizontal $4T$ relations".

The next theorem is technical, but provides the foundation required
to prove the main results here; the proof is given in \fullref{sec:Proof of Theorem.PSigma.n.fibrations}.

\begin{thm}\label{thm:PSigma.n.fibrations} There exist homomorphisms
$$\pi\co P\varSigma_n \to\ P\varSigma_{n-1}$$ defined by
\[
\pi(\chi_{k,i}) =
\begin{cases}
\chi_{k,i} & \text{if $i < n$, and $k < n$,}\\
1 & \text{if $i = n$ or $k=n$.}
\end{cases}
\] The homomorphism $\pi\co  P\varSigma_n \to P\varSigma_{n-1}$ is an epimorphism.
The kernel of $\pi$ denoted $K_n$ is generated by the elements
$\chi_{n,i}$ and $\chi_{j,n}$ for $1 \leq i,j \leq n-1$.
Furthermore, this extension is split and the conjugation action of
$P\Sigma_{n-1}$ on $H_1(K_n)$ is trivial.

In addition, the homomorphism $\pi\co P\varSigma_n \to\
P\varSigma_{n-1}$ restricts to a homomorphism
$$\pi|_{P\varSigma_{n}^+}\co P\varSigma_n^+ \to\ P\varSigma_{n-1}^+$$ which is an epimorphism. The kernel of
$\pi|_{P\varSigma_n^+}$ denoted $K_n^+$ is a free group generated by
the elements $\chi_{n,i}$ for $1 \leq i \leq n-1$. Furthermore, this
extension is split and the conjugation action of
$P\varSigma_{n-1}^+$ on $H_1(K_n^+)$ is trivial.
\end{thm}

{\bf Remark 1}\qua After this paper was submitted, the authors learned
that the result in \fullref{thm:PSigma.n.fibrations} stating that
$\pi\co P\varSigma_n \to\ P\varSigma_{n-1}$ is a split epimorphism was
proved earlier by Bardakov in \cite{bardakov} where other natural
properties are developed. In addition, the feature that $K_n$ is a
semi-direct product as stated in \fullref{thm:PSigma.n.Lie.algebras}
below was also proved in \cite{bardakov}.

Provided that it is clear from the context, the notation
$\chi_{k,i}$ is used ambiguously to denote both the element
$\chi_{k,i}$ in $P\varSigma_n$, or in $P\varSigma_n^+$ when defined,
as well as the equivalence class of $\chi_{k,i}$ in
$gr^1(P\varSigma_n)= H_1(P\varSigma_n)$ or in $H_1(P\varSigma_n^+)$
when defined. Partial information concerning the Lie algebra
$gr^*(P\varSigma_n)$ is given next.

\begin{thm}\label{thm:PSigma.n.Lie.algebras} There is a
split short exact sequence of Lie algebras
$$0 \to gr^*(K_n) \to  gr^*(P\varSigma_n) \to  gr^*(P\varSigma_{n-1}) \to
0.$$ The relations (1)--(4) are satisfied on the level of Lie
algebras:
\begin{enumerate}
\item If $\{i,j\} \cap \{s,t\} = \emptyset$, then
  $[\chi_{j,i}, \chi_{s,t}] = 0$.
  \item If $ i, j, k$ are distinct, then $[\chi_{i,k},\chi_{i,j}+ \chi_{k,j}] = 0.$
\item If $ i,j,k$ are distinct, the element $[\chi_{k,i}, \chi_{j,i} + \chi_{j,k}]$
    is non-zero.
\item If $ i,j,k$ are distinct, $[ \chi_{i,j}, \chi_{k,j}] = 0$.
\end{enumerate}

Furthermore, there is a split epimorphism $$\gamma\co  K_n \to
\oplus_{n-1}\mathbb Z $$ with kernel denoted $\Lambda_n$ together
with a split short exact sequence of Lie algebras
\[
\begin{CD}
0 @>{}>> \mathbb L_n @>{}>> gr^*(K_n) @>{gr^*(\gamma)}>>
\oplus_{n-1}\mathbb Z @>{}>> 0
\end{CD}
\] where $\mathbb L_n$ is the Lie algebra kernel of $gr^*(\gamma)$.
\end{thm}%\vskip .04in

The same methods give a complete description for the cohomology
algebra of $P\varSigma_n^+$ as well as the Lie algebra
$gr^*(P\varSigma_n^+)$. A further application to be given later is a
substantial contribution to the cohomology of each of the Johnson
filtrations of $IA_n$. To express the answers, the notation
$\chi_{k,i}^*$ is used to denote the dual basis element to
$\chi_{k,i}$, namely
\[
\chi_{k,i}^*(\chi_{s,t}) =
\begin{cases}
1
  &   \text{if $k=s$ and $i=t$ } \\
0
  & \text{otherwise. } \\
\end{cases}
\]

\begin{thm}\label{thm:PSigma_{n,+}}
The cohomology algebra of $P\varSigma_n^+$ satisfies the following
  properties.
            \begin{enumerate}
            \item Each graded piece $H^k(P^+\Sigma_n)$ is a finitely
generated, torsion-free abelian group.
            \item If $ 1 \leq k \leq n$, a basis for $H^kP\varSigma_n^+$ is given by
            $\chi_{i_1,j_1}^*\cdot \chi_{i_2,j_2}^*\cdots \chi_{i_k,j_k}^*$ where $2
            \leq i_1 < i_2 < \cdots < i_k \leq n,$ and $1\leq j_t < i_t $ 
for all $t$.
            \item A complete set of relations (assuming graded commutativity, 
and associativity)
            is given by
        \begin{itemize}
                \item ${\chi_{i,k}^*}^2 = 0$ for all $i > k$, and
                \item $\chi_{i,j}^*[\chi_{i,k}^*- \chi_{j,k}^*]=0$ for $k < j < i$.
        \end{itemize}
            \end{enumerate}
The Lie algebra obtained from the descending central series for
$P\varSigma_n^+$, $gr^*(P\varSigma_n^+)$, is additively isomorphic
to a direct sum of free sub-Lie algebras $$\bigoplus_{2 \leq k \leq n}
L[\chi_{k,1}, \chi_{k,2}, \cdots, \chi_{k,k-1}],$$ with
\begin{itemize}
    \item $[ \chi_{k,j}, \chi_{s,t}] = 0$ if $\{j,k\} \cap \{s,t \}=
\emptyset$,
    \item $[ \chi_{k,j}, \chi_{s,j}] = 0$ if $\{s,k\} \cap \{j \}= \emptyset$ and

    \item $[\chi_{i,k},\chi_{i,j}+ \chi_{k,j}] = 0$ for $ j < k < i.$
\end{itemize}
\end{thm}

{\bf Remark 2}\qua The structure of the cohomology algebra described in
Theorem~1.4 corresponds to the structure of the algebra $H^*(P\varSigma_n,
\Z)$, as given in Jensen, McCammond and Meier \cite{jmm}, under the
map induced by the canonical inclusion $P\varSigma_n^+ \subset
P\varSigma_n$.  Namely, the elements $\chi_{k,i}^*$ with $k < i$ are
mapped to zero and relations of $H^*(P\varSigma_n, \Z)$ become the
relation (3) of Theorem~1.4.

{\bf Remark 3}\qua \fullref{thm:PSigma_{n,+}} does not rule out the
possibility that $P\varSigma_n^+$ is isomorphic to the pure braid
group $P_n$. Notice that $P\varSigma_3^+$ is isomorphic to $\mathbb
Z \times F[\chi_{3,1},\chi_{3,2}]$ and thus $P_3$ where a generator
of $\mathbb Z$ is given by $\chi_{2,1}\cdot \chi_{3,1}$. In fact,
\fullref{thm:PSigma_{n,+}} implies that after suspending the
classifying spaces of both $P\varSigma_n^+$ and $P_n$ exactly once,
these suspended classifying spaces are homotopy equivalent
as they are both finite bouquets of spheres
with the same dimensions.

Let $U[\mathcal L]$ denote the universal enveloping algebra of a Lie
algebra $\mathcal L$. Since the Euler--Poincar\'e series for
$U[L[\chi_{k,1}, \chi_{k,2}, \cdots, \chi_{k,k-1}]]$ is
$1/(1-(k-1)t)$, the next corollary follows at once.

\begin{cor}\label{cor:PSigma_{n,+}}
The Euler--Poincar\'e series for $U[P\varSigma_n^+]$ is equal to
$$\prod_{1 \leq k \leq n-1} 1/(1-kt).$$
\end{cor}

\begin{defin}
Let $\mathbb M_n$  denote the smallest subalgebra of $H^*(IA_n;
\mathbb Z)$ such that
\begin{enumerate}
\item $\mathbb M_n$ surjects to $H^*(P\varSigma_n^+)$ (that such a surjection
exists follows from \fullref{thm:A(n)}) and
\item $\mathbb M_n$ is closed respect to the conjugation action of $GL(n, \mathbb Z)$ on
$H^*(IA_n; \mathbb Z)$.
\end{enumerate}
\end{defin}

\begin{thm}\label{thm:A(n)}
The natural inclusion $j\co P\varSigma_n^+ \to IA_n$ composed with the
abelianization map $A\co  IA_n \to IA_n/[IA_n,IA_n] =
\oplus_{n\binom{n}{2}} \mathbb Z $ given by
\[
\begin{CD}
P\varSigma_n^+   @>{j}>> IA_n  @>{A}>>  IA_n/[IA_n,IA_n]=H_1(IA_n)
\end{CD}
\] induces a split epimorphism in integral cohomology. Thus the
integral cohomology of $P\varSigma_n^+$ is a direct summand of the
integer cohomology of $IA_n$ (a summand which is not invariant
under the action of $GL(n, \mathbb Z)$). Furthermore, the image of
$$A^*\co H^*(IA_n/[IA_n,IA_n]) \to H^*(IA_n)$$ contains $\mathbb M_n$.

In addition, the suspension of the classifying space
$BP\varSigma_n^+ $, $\Sigma(BP\varSigma_n^+ )$, is a retract of
$\Sigma(BIA_n)$ and there is an induced map $$ \theta\co  BIA_n \to
\Omega \Sigma (BP\varSigma_n^+ )$$ which factors the Freudenthal
suspension $ E\co  BP\varSigma_n^+  \to \Omega \Sigma(BP\varSigma_n^+
)$ given by the composite $$BP\varSigma_n^+ \to  BIA_n \to \Omega
\Sigma (BP\varSigma_n^+).$$
\end{thm}

{\bf Remark~4}\qua Properties of the image of $A^*$ are addressed in
work of N~Kawazumi \cite{ka}, T~Sakasai \cite{s}, T~Satoh
\cite{sa} and A~Pettet \cite{pe}. An analogous map  $BIA_n \to
\Omega \Sigma (BP\varSigma_n)$ is constructed in work of C~Jensen,
J~McCammond and J~Meier \cite{jmm}.

%%%%%%%%%%%%%%%%%%%%%%%%%%%%%%%%%%%%%%%%%%%%
%%ADDED 20 MAY 2006
%%%%%%%%%%%%%%%%%%%%%%%%%%%%%%%%%%%%%%%%%%%

The remainder of this introduction is devoted to the structure of
the {\it braid--permutation group} $\mathrm{BP}_n$ introduced by
R~Fenn, R~Rim\'anyi and C~Rourke \cite{FRR2}. The group
$\mathrm{BP}_n$ is defined as the subgroup of $\Aut(F_n)$ generated
by $\xi_i$ and $\sigma_i$, where the action of an element $\phi$ in
$\Aut(F_n)$ is from the right with
$$
(x_j)\xi_i=\left\{
\begin{array}{cc}
x_{i+1}&j=i,\\
x_i&j=i+1,\\
x_j&{\rm otherwise;}
\end{array}
\right.
$$
$$
(x_j)\sigma_i=\left\{
\begin{array}{cc}
x_{i+1}&j=i,\\
x_{i+1}^{-1}x_ix_{i+1}&j=i+1,\\
x_j&{\rm otherwise.}
\end{array}
\right.
$$
The group $\mathrm{BP}_n$ is presented by the set of generators
$\xi_i$ and $\sigma_i$ for $1\leq i\leq n-1$, and by  the relations:
\begin{equation}\label{relation1}
\left\{
\begin{array}{cc}
\xi_i^2=1,& \\
\xi_i\xi_j=\xi_j\xi_i& |i-j|>1,\\
\xi_i\xi_{i+1}\xi_i=\xi_{i+1}\xi_i\xi_{i+1};
\end{array}
\right.
\end{equation}
\begin{equation}\label{relation2}
\left\{
\begin{array}{cc}
\sigma_i\sigma_j=\sigma_j\sigma_i& |i-j|>1,\\
\sigma_i\sigma_{i+1}\sigma_i=\sigma_{i+1}\sigma_i\sigma_{i+1};
\end{array}
\right.
\end{equation}
\begin{equation}\label{relation3}
\left\{
\begin{array}{cc}
\xi_i\sigma_j=\sigma_j\xi_i& |i-j|>1,\\
\xi_i\xi_{i+1}\sigma_i=\sigma_{i+1}\xi_i\xi_{i+1},\\
\sigma_i\sigma_{i+1}\xi_i=\xi_{i+1}\sigma_i\sigma_{i+1},\\
\end{array}
\right.
\end{equation}

\smallskip

The group $\mathrm{BP}_n$ is also characterized \cite{FRR2} as the
the subgroup of $\Aut(F_n)$ consisting of automorphism $\phi\in
\Aut(F_n)$ of {\it permutation--conjugacy type} which satisfy
\begin{equation}
%%\label{b_p_type}
(x_i)\phi=w_i^{-1}x_{\lambda(i)}w_i
\end{equation}
for some word $w_i\in F_n$ and permutation $\lambda \in \Sigma_n$
the symmetric group on $n$ letters.

\begin{thm}\label{thm:braid.permutation}  The group $\mathrm{BP}_n$
is the semi-direct product of the symmetric group on $n$--letters
$\Sigma_n$ and the group $P\varSigma_n$ with a split extension
\[
\begin{CD}
1 @>{}>>  P\varSigma_n @>{}>> \mathrm{BP}_n  @>{}>> \Sigma_n @>{}>>
1.
\end{CD}
\]
\end{thm}

Theorems \ref{thm:PSigma.n.Lie.algebras}, \ref{thm:PSigma_{n,+}} and
\ref{thm:braid.permutation} provide some information about the
cohomology as well as Lie algebras associated to $\mathrm{BP}_n$.

The authors take this opportunity to thank Toshitake Kohno, Nariya
Kawazumi,\break Shigeyuki Morita, Dai Tamaki as well as other friends for
this very enjoyable opportunity to participate in this conference.
The authors would like to thank Benson Farb for his interest in this
problem. The authors would also like to thank Allen Hatcher for his
careful reading and interest in this article.

The first author is especially grateful for this mathematical
opportunity to see friends as well as to learn and to work on
mathematics with them at this conference.

\section{Projection maps $P\varSigma_n \to P\varSigma_{n-1}$}

Consider the map $$p\co F_n \to\ F_{n-1}$$ defined by the formula
\[
p(x_j)=
\begin{cases}
x_j & \text{if $j \leq n-1$,}\\
1 & \text{if $j = n$.}
\end{cases}
\] In addition, let $\nu\co P\varSigma_n^+ \to P\varSigma_n$ denote the 
natural inclusion.

\begin{thm} \label{thm:restrictions}
The projection maps $p\co F_n \to\ F_{n-1}$ induce homomorphisms
$$\pi\co P\varSigma_n \to\ P\varSigma_{n-1}$$ given by
\[
\pi(\chi_{k,i}) =
\begin{cases}
\chi_{k,i} & \text{if $i < n$, and $k < n$,}\\
1 & \text{if $i = n$ or $k=n$.}
\end{cases}
\] Furthermore, these homomorphisms restrict to
$$\pi\co P\varSigma_n^+ \to\ P\varSigma_{n-1}^+$$ together with a commutative diagram:
\[
\begin{CD}
P\varSigma_n^+ @>{\pi}>> P\varSigma_{n-1}^+    \\
  @V{\nu}VV                    @VV{\nu}V \\
P\varSigma_n @>{\pi}>> P\varSigma_{n-1}
\end{CD}
\]
\end{thm}
\begin{proof}
Consider the following commutative diagrams:
\[
\begin{CD}
F_n @>{\chi_{n,j}}>> F_n    \\
  @V{p}VV                    @VV{p}V \\
F_{n-1}  @>{1}>> F_{n-1}
\end{CD}
\]
\[
\begin{CD}
F_n @>{\chi_{j,n}}>> F_n    \\
  @V{p}VV                    @VV{p}V \\
F_{n-1}  @>{1}>> F_{n-1}
\end{CD}
\]
If $k,j< n$, then the following diagram commutes:
\[
\begin{CD}
F_n @>{\chi_{j,k}}>> F_n    \\
  @V{p}VV                    @VV{p}V \\
F_{n-1}  @>{\chi_{j,k}}>> F_{n-1}
\end{CD}
\]
Thus, if any of $$i,k <n \qua \hbox{or} \qua
i= n\qua \hbox{or} \qua k = n$$ hold, then
the functions $\chi_{k,i}$ restrict to isomorphisms of $F_{n-1}$. The
restriction is evidently compatible with composition of
isomorphisms. Hence there is an induced homomorphism
$$\pi\co P\varSigma_n \to\ P\varSigma_{n-1}$$ given by
\[
\pi(\chi_{k,i}) =
\begin{cases}
\chi_{k,i} & \text{if $i < n$ and $k < n$,}\\
1 & \text{if $i = n$ or $k=n$.}
\end{cases}
\]
These homomorphisms are compatible with the inclusion maps
$\nu\co P\varSigma_n^+ \to P\varSigma_n$, and the theorem follows.
\end{proof}

\section{On automorphisms of $P\varSigma_n$, and $P\varSigma_n^+$}

The conjugation action of $\Aut(F_n)$ on itself restricted to certain
natural subgroups of $\Aut(F_n)$ has $P\varSigma_n$, and
$P\varSigma_n^+$ as characteristic subgroups. The purpose of this
section is to give two such natural subgroups. One of these
subgroups is used below to determine the relations in the cohomology
algebra for $P\varSigma_n^+$.

A choice of generators for $\Aut(F_n)$ is listed next. Let $\sigma$
denote an element in the symmetric group on $n$ letters $\Sigma_n$
which acts naturally on $F_n$ by permutation of coordinates with
$\xi_i$ given by the transposition $(i,i+1)$ (see the end of the previous 
section). Thus
\[
\xi_i(x_j)=
\begin{cases}
x_j & \text{if $\{j\} \cap \{i,i+1\} = \emptyset$},\\
x_{i+1} & \text{ if $j = i$ and} \\
x_i & \text{if $j = i+1$.}
\end{cases}
\]
Next consider $$\tau_i\co  F_n \to F_n$$ which sends $x_i$ to
$x_i^{-1}$, and fixes $x_j$ for $j \neq i$. Thus
\[
\tau_i(x_j)=
\begin{cases}
x_j & \text{if $j \neq i $ and}  \\
x_i^{-1} & \text{if $j = i$.}
\end{cases}
\]
The elements $\tau_i$, and $\xi_j$ for $ 1 \leq i,j \leq n$
generate the ``signed permutation group" in $\Aut(F_n)$; this group,
also known as the wreath product $\Sigma_n \wr \mathbb Z/ 2 \mathbb
Z$, embeds in $GL(n,\mathbb Z)$ via the natural map $\Aut(F_n) \to
GL(n,\mathbb Z)$. In this wreath product, a Coxeter group of type
$B_n$, the elements $\tau_i$ can be expressed in terms of $\tau_1$
and $\xi_i$.

Let $\delta$ denote the automorphism of $F_n$ which sends $x_1$ to
$x_1 x_2$ while fixing $x_i$ for $i > 1$. It follows from \cite{mks}
that $\Aut(F_n)$ is generated by the elements
\begin{itemize}
    \item $\xi_i$ for $1 \leq i \leq n$,
    \item $\tau_i$ for $1 \leq i \leq n$ and
    \item $\delta$.
\end{itemize}
It is natural to consider the action of some of these elements on
$P\varSigma_n$, and $P\varSigma_n^+$ by conjugation.
\begin{prop} \label{prop:automorphisms}
Subgroups of the automorphism groups of $P\varSigma_n^+$, and
$P\varSigma_n$ are listed as follows.
\begin{enumerate}
%%OK
\item Conjugation by the elements $\tau_i$ for $1 \leq i \leq n$
leaves $P\varSigma_n^+$ invariant. Thus $\oplus_n \mathbb Z/ 2
\mathbb Z$ is isomorphic to a subgroup of $\Aut(P\Sigma_{n}^+)$ with
induced monomorphisms $$\theta\co \oplus_n \mathbb Z/ 2 \mathbb Z \to
\Aut(P\varSigma_{n}^+)$$ obtained by conjugating an element in
$P\varSigma_{n}^+$ by the $\tau_i$. The conjugation action of
$\tau_i$ on the elements $\chi_{s,t}$ is specified by the formulas
    \[
    \tau_i(\chi_{s,t})\tau_i^{-1}=
    \begin{cases}
    \chi_{s,t} & \text{if $\{i\} \cap \{s,t\} = \emptyset$},\\
    \chi_{s,i}^{-1} & \text{ if $t = i$},\\
    \chi_{i,t} & \text{if $s = i$}
    \end{cases}
    \] for $s > t$.

%%OK
\item Conjugation by the elements $\tau_i$ and $\xi_i$ for $1 \leq i \leq n$
leaves $P\varSigma_n$ invariant. Thus $\Sigma_n \wr \mathbb Z/
2\mathbb Z$ is isomorphic to a subgroup of $\Aut(P\Sigma_{n})$ with
induced monomorphisms $$\theta\co \Sigma_n \wr \mathbb Z/ 2\mathbb Z
\to \Aut(P\varSigma_n)$$ obtained by conjugating an element in
$P\varSigma_{n}$ by the $\tau_i$ and $\xi_i$. The conjugation
action of $\tau_i$ and $\xi_i$ on the $\chi_{s,t}$ is specified
by the formulas
    %%OK
    \[
    \tau_i(\chi_{s,t})\tau_i^{-1}=
    \begin{cases}
    \chi_{s,t} & \text{if $\{i\} \cap \{s,t\} = \emptyset$},\\
    \chi_{s,i}^{-1} & \text{ if $t = i$},\\
    \chi_{i,t} & \text{if $s = i$,}
    \end{cases}
    \] and
%%OK
    \[
    \xi_i(\chi_{s,t})\xi_i^{-1}=
    \begin{cases}
    \chi_{s,t} & \text{if $\{i,i+1\} \cap \{s,t\} = \emptyset$ },\\
    \chi_{i+1,t} & \text{ if $s = i$ and $t \neq i+1$},\\
    \chi_{i+1,i} & \text{ if $s = i$ and $t = i+1$},\\
    \chi_{i,t} & \text{ if $s = i+1$ and $t \neq i$},\\
    \chi_{i,i+1} & \text{ if $s = i+1$ and $t = i$},\\
    \chi_{s,i+1} & \text{ if $s \neq i+1$ and $t = i$}, \\
    \chi_{s,i} & \text{ if $s \neq i$ and $t = i+1$},
    \end{cases}
    \] for all $s \neq t$.
\end{enumerate}
\end{prop}

\begin{proof}
The proof of this proposition is a direct computation with details
omitted.
\end{proof}

\section{On certain subgroups of $IA_n$}
\label{sec:On certain subgroups of $IA_n$}

The purpose of this section is to consider the subgroups

\begin{enumerate}
  \item  $\mathcal G_n$ of $P\varSigma_{n}$ generated by
the elements $\chi_{n,i}$ and $\chi_{j,n}$ for $1 \leq i, j \leq
n-1$ and
  \item $\mathcal G_n^+$ of $P\varSigma_{n}^+$ generated by
the elements $\chi_{n,i}$ for $1 \leq i\leq n-1$.
\end{enumerate}

\begin{prop} \label{prop:normal.subgroups}
\begin{enumerate}
\item The following relations are satisfied in $P\Sigma_n.$

\begin{enumerate}
\item[(i)]  $\chi_{i,j}^{-1}\cdot \chi_{n,j}\cdot \chi_{i,j}=
  \chi_{n,j}$ with $i,j < n$.

   \item[(ii)] $\chi_{i,k}^{-1}\cdot \chi_{n,j}\cdot \chi_{i,k} =
  \chi_{n,j}$ with $\{i,k\} \cap \{n,j\}= \emptyset$.

 \item[(iii)] $\chi_{j,k}^{-1}\cdot \chi_{n,j}\cdot \chi_{j,k} =
  \chi_{n,k} \cdot \chi_{n,j} \cdot \chi_{n,k}^{-1}$  with $k,j <n$.
  \end{enumerate}

  \medskip

\begin{enumerate}

  \item[(iv)] $\chi_{i,k}^{-1} \cdot \chi_{j,n} \cdot \chi_{i,k} = \chi_{j,n} $ with $\{i,k\} \cap \{n,j \}= \emptyset$.

\item[(v)] $\chi_{j,i}^{-1}\cdot \chi_{j,n}\cdot \chi_{j,i}=
\chi_{n,i}\cdot \chi_{j,n}\cdot \chi_{n,i}^{-1}$  with $i,j<n$.

  \item[(vi)] $\chi_{i,j}^{-1} \cdot \chi_{j,n}\cdot \chi_{i,j} = (\chi_{n,j} \cdot
\chi_{i,n}^{-1} \cdot \chi_{n,j}^{-1}) \cdot (\chi_{i,n} \cdot
\chi_{j,n})$ with $i,j < n$.
\end{enumerate}

  \item The group $\mathcal G_n$ is a normal subgroup of $P\varSigma_{n}$
and is the kernel of the projection $$\pi\co P\Sigma_n \to\
P\Sigma_{n-1}.$$ Thus $\mathcal G_n = K_n$ as given in \fullref{thm:restrictions}.

  \item
The group $\mathcal G_n^+$ is a normal subgroup of
$P\varSigma_{n}^+$ and is the kernel of the projection
$$\pi|_{P\varSigma_{n}^+}\co P\Sigma_n^+ \to\ P\Sigma_{n-1}^+.$$
Thus $\mathcal G_n^+ = K_n^+$ as given in \fullref{thm:restrictions}.

\end{enumerate}
\end{prop}

\begin{proof}

Consider the elements $\chi_{i,j}^{-1}\cdot \chi_{t,n}\cdot
\chi_{i,j}$ and $\chi_{i,j}^{-1}\cdot \chi_{n,t}\cdot \chi_{i,j}$
for various values of $t$ together with McCool's relations as listed
in \fullref{thm:McCool relations}:
\begin{enumerate}
\item $\chi_{i,j}\cdot \chi_{k,j}\cdot \chi_{i,k} = \chi_{i,k} \cdot
\chi_{i,j} \cdot \chi_{k,j}$ for $i,j,k$ distinct.
\item $[\chi_{k,j}, \chi_{s,t}] = 1$ if $\{j,k\} \cap \{s,t \}=
\emptyset$.
\item $[ \chi_{i,j},\chi_{k,j}] = 1$ for
$i,j,k$ distinct.
\end{enumerate}

The verification of part $(1)$ of the proposition breaks apart into
six natural cases where the first three are given by the conjugation
action on $\chi_{n,j}$ while the second three are given by the
conjugation action on $\chi_{j,n}$.

\begin{enumerate}
  %%2,2 OK
\item[(i)]  $\chi_{i,j}^{-1}\cdot \chi_{n,j}\cdot \chi_{i,j}=
  \chi_{n,j}$ by formula (3) with $i,j < n$.
%%empty OK
   \item[(ii)] $\chi_{i,k}^{-1}\cdot \chi_{n,j}\cdot \chi_{i,k} =
  \chi_{n,j}$ by formula (2) with $\{i,k\} \cap \{n,j\}= \emptyset$.
%%2,1 OK
 \item[(iii)] $\chi_{j,k}^{-1}\cdot \chi_{n,j}\cdot \chi_{j,k} =
  \chi_{n,k} \cdot \chi_{n,j} \cdot \chi_{n,k}^{-1}$ by formulas (1) and (3) with $k,j <n$.
  \end{enumerate}

  \medskip

%%%%%%%%%%%%%%%%%%%%%%%%%%%%%%%%%%%%%%%%%%%%%%%
\begin{enumerate}
  %%empty OK
  \item[(iv)] $\chi_{i,k}^{-1} \cdot \chi_{j,n} \cdot \chi_{i,k} = \chi_{j,n} $ by
formula (2) with $\{i,k\} \cap \{n,j \}= \emptyset$.

%% 1-1 OK
\item[(v)] $\chi_{j,i}^{-1}\cdot \chi_{j,n}\cdot \chi_{j,i}=
\chi_{n,i}\cdot \chi_{j,n}\cdot \chi_{n,i}^{-1}$ by formulas (1) and
(3) with $i,j<n$.

  %%2-1 NOT OK
  \item[(vi)] $\chi_{i,j}^{-1} \cdot \chi_{j,n}\cdot \chi_{i,j} = (\chi_{n,j} \cdot
\chi_{i,n}^{-1} \cdot \chi_{n,j}^{-1}) \cdot (\chi_{i,n} \cdot
\chi_{j,n})$ by formulas (1) and (3) with $i,j < n$.
\end{enumerate}

A sketch of formula (iv) is listed next for convenience of the
reader. Assume that $i,j,n$ are distinct.
\begin{enumerate}
  \item[(a)] $\chi_{i,n}\cdot \chi_{j,n}\cdot \chi_{i,j} = \chi_{i,j} \cdot
\chi_{i,n} \cdot \chi_{j,n}$.
  \item[(b)] $\chi_{j,n}\cdot \chi_{i,j} = \chi_{i,n}^{-1}\cdot \chi_{i,j} \cdot
\chi_{i,n} \cdot \chi_{j,n}$.
  \item[(c)] $\chi_{i,j}^{-1} \cdot \chi_{j,n}\cdot \chi_{i,j} = 
(\chi_{i,j}^{-1} \cdot \chi_{i,n}^{-1}\cdot \chi_{i,j}) \cdot
(\chi_{i,n} \cdot \chi_{j,n})$.
\item[(d)] By formula (v) above, 
$\chi_{i,j}^{-1}\cdot \chi_{i,n}\cdot \chi_{i,j}= \chi_{n,j}\cdot
\chi_{i,n}\cdot \chi_{n,j}^{-1}$.
\item[(e)] Thus $\chi_{i,j}^{-1} \cdot \chi_{j,n}\cdot \chi_{i,j} = (\chi_{n,j}\cdot
\chi_{i,n}\cdot \chi_{n,j}^{-1}) \cdot (\chi_{i,n} \cdot
\chi_{j,n})$ and formula (vi) follows.
\end{enumerate}

Thus \begin{enumerate}
       \item $\mathcal G_n$ is a normal subgroup of
$P\varSigma_{n}$ and
       \item $\mathcal G_n^+$ is a normal subgroup of
$P\varSigma_{n}^+$
     \end{enumerate}
by inspection of the previous relations (i--vi).

Next notice that $\mathcal G_n$ is in the kernel of the projection
$\pi\co P\Sigma_n \to\ P\Sigma_{n-1}$. 
Denote by $\sigma$ the canonical section 
$\sigma\co P\varSigma_{n-1} \to\ P\varSigma_{n}$ with
$\sigma(\chi_{j,i}) = \chi_{j,i}$. 
Every element $W$
in $P\Sigma_n$ is equal to a product $W_{n-1} \cdot X_n$ where
$W_{n-1}$ is in the image of the section $\sigma$ applied to
$P\Sigma_{n-1}$ and $X_n$ is in $\mathcal G_n$ by inspection of the
relations relations (i--vi). Thus, the kernel of $\pi$ is generated
by all conjugates of the elements $\chi_{n,i}$ and $\chi_{j,n}$ for
$1 \leq i, j \leq n-1$ which coincides with $\mathcal G_n$. A
similar assertion and proof applies to $\mathcal G_n^+$ using
relations (iv, v, vi).
\end{proof}

The following result \cite{k,fr} will be used below to derive the
structure of certain Lie algebras in this article.
\begin{thm} \label{thm:exact sequences of Lie algebras}
Let $$1\to A \to B \to C \to 1$$ be a split short exact sequence of
groups for which conjugation by $C$ induces the trivial action on
$H_1(A)$. Then there is a split short exact sequence of Lie algebras
$$0 \to\ gr^*(A) \to\ gr^*(B) \to\ gr^*(C) \to\ 0.$$
\end{thm}

To apply \fullref{thm:exact sequences of Lie algebras}, features
of the local coefficient system in homology for the projection maps
$\pi\co P\Sigma_n \to\ P\Sigma_{n-1}$ and
$\pi|_{P\varSigma_{n}^+}\co P\Sigma_n^+ \to\ P\Sigma_{n-1}^+$ are
obtained next.

\begin{prop} \label{prop:local.coefficients}
\begin{enumerate}
\item The natural conjugation action of $P\Sigma_{n-1}$ on $H_1(K_n)$ is trivial.
Thus there is a split short exact sequence of Lie algebras
$$0 \to\ gr^*(K_n) \to\  gr^*(P\Sigma_n) \to\  gr^*(P\Sigma_{n-1}) \to\ 0.$$

\item The natural conjugation action of $P\Sigma_{n-1}^+$ on $H_1(K_n^+)$ is trivial.
Thus there is a split short exact sequence of Lie algebras
$$0 \to\ gr^*(K_n^+) \to\  gr^*(P\Sigma_n^+) \to\  gr^*(P\Sigma_{n-1}^+) \to\ 0.$$
\end{enumerate}
\end{prop}

\begin{proof}
As before, consider the elements $\chi_{i,j}^{-1}\cdot
\chi_{n,t}\cdot \chi_{i,j}$ and $\chi_{i,j}^{-1}\cdot
\chi_{t,n}\cdot \chi_{i,j}$ together with McCool's relations as
given in \fullref{prop:normal.subgroups} to obtain the
following formulas.

\begin{enumerate}
\item[(i)]  $\chi_{i,j}^{-1}\cdot \chi_{n,j}\cdot \chi_{i,j} = \chi_{n,j}$ 
with $i,j < n$.

\item[(ii)] $\chi_{i,k}^{-1}\cdot \chi_{n,j}\cdot \chi_{i,k} = \chi_{n,j}$ 
with $\{i,k\} \cap \{n,j\}= \emptyset$.

\item[(iii)] $\chi_{j,k}^{-1}\cdot \chi_{n,j}\cdot \chi_{j,k} = \chi_{n,j}\cdot 
(\chi_{n,j}^{-1}\cdot \chi_{n,k}
\cdot \chi_{n,j} \cdot \chi_{n,k}^{-1})$  with $k,j <n$. Thus
$\chi_{j,k}^{-1}\cdot \chi_{n,j}\cdot \chi_{j,k} = \chi_{n,j}\cdot
([\chi_{n,j}^{-1}, \chi_{n,k}]).$

\end{enumerate}

\medskip

\begin{enumerate}
\item[(iv)] $\chi_{i,k}^{-1} \cdot \chi_{j,n} \cdot \chi_{i,k} = \chi_{j,n} $ 
with $\{i,k\} \cap \{n,j \}= \emptyset$.
\item[(v)] $\chi_{j,i}^{-1}\cdot \chi_{j,n}\cdot \chi_{j,i}=
\chi_{n,i}\cdot \chi_{j,n}\cdot \chi_{n,i}^{-1} = \chi_{j,n}\cdot
 \chi_{j,n}^{-1}\cdot  \chi_{n,i}\cdot \chi_{j,n}\cdot \chi_{n,i}^{-1}  $  
with $i,j<n$.
 Thus $\chi_{j,i}^{-1}\cdot \chi_{j,n}\cdot \chi_{j,i}=
 \chi_{j,n}\cdot [\chi_{j,n}^{-1},\chi_{n,i}]$.
\item[(vi)] $\chi_{i,j}^{-1} \cdot \chi_{j,n}\cdot \chi_{i,j} = (\chi_{n,j} \cdot
\chi_{i,n}^{-1} \cdot \chi_{n,j}^{-1}) \cdot (\chi_{i,n} \cdot
\chi_{j,n})= [\chi_{n,j},\chi_{i,n}^{-1}]\cdot \chi_{j,n}$ with $i,j
< n$.
\end{enumerate}

It then follows that the conjugation action of $\chi_{s,t}$ on the
class of either $\chi_{n,j}$ or $\chi_{j,n}$ in $H_1(K_n)$ fixes
that class (in $H_1(K_n)$). Thus there is a short exact sequence of
Lie algebras $0 \to\ gr^*(K_n) \to\ gr^*(P\Sigma_n) \to\
gr^*(P\Sigma_{n-1}) \to\ 0$ by \fullref{thm:exact sequences of
Lie algebras}

A similar assertion and proof follows for $H_1(K_n^+)$ by inspection
of formulas (iv,v,vi). The proposition follows.
\end{proof}

\section[Proof of \ref{thm:PSigma.n.fibrations}]{Proof of \fullref{thm:PSigma.n.fibrations}}
\label{sec:Proof of Theorem.PSigma.n.fibrations}

The first part of \fullref{thm:PSigma.n.fibrations} that the
projection map $\pi\co P\varSigma_n \to\ P\varSigma_{n-1}$ is an
epimorphism with kernel generated by $\chi_{n,i}$ and $\chi_{j,n}$
for $1 \leq i,j \leq n-1$ follows from \fullref{prop:normal.subgroups}. Furthermore, this extension is split by
the section $\sigma\co P\varSigma_{n-1} \to\ P\varSigma_{n}$. That the 
local coefficient system
%% added 19 october 2006
is trivial for $H_1(K_n)$ follows from \fullref{prop:local.coefficients}. (Note that it is possible that the
local coefficient system is non-trivial for higher dimensional
homology groups of $K_n$.) Similar properties are satisfied for
$\pi|_{P\varSigma_{n}^+}\co P\varSigma_n^+ \to\ P\varSigma_{n-1}^+$ by
Propositions \ref{prop:normal.subgroups} and
\ref{prop:local.coefficients}.

Consider the free group $F_{n-1}$ with generators $x_1, \cdots,
x_{n-1}$. There is a homomorphism $$\Phi_{n-1}\co F_{n-1}\to K_n^+$$
obtained by defining $$\Phi_{n-1}(x_i) =\chi_{n,i}.$$ This
homomorphism is evidently a surjection. To check that the subgroup
$K_n^+$ is a free group generated by the elements $\chi_{n,i}$ for
$1 \leq i \leq n-1$, it suffices to check that $\Phi_{n-1}$ is a
monomorphism.

Observe that if $$W= x_{i_1}^{\epsilon_1}\cdot x_{i_2}^{\epsilon_2}
\cdots x_{i_r}^{\epsilon_r}$$ for $\epsilon_i = \pm 1$ is a word in
$F_{n-1}$, then $\Phi_{n-1}(W)$ is an automorphism of $F_n$ with
$$\Phi_{n-1}(W)(x_n) = W\cdot x_n \cdot W^{-1}.$$ Thus if
$W$ is in the kernel of $\Phi_{n-1}$, then $W\cdot x_n \cdot
W^{-1}=x_n$. Furthermore, if $W$ is in $F_{n-1}$, then $W = 1$.
\fullref{thm:PSigma.n.fibrations} follows.

\section[Proof of \ref{thm:PSigma.n.Lie.algebras}]{Proof of \fullref{thm:PSigma.n.Lie.algebras}}
\label{sec:Proof of Theorem.PSigma.n.Lie.algebras}

\fullref{thm:PSigma.n.Lie.algebras} states that there a split
short exact sequence of Lie algebras $$0 \to gr^*(K_n) \to
gr^*(P\varSigma_n) \to gr^*(P\varSigma_{n-1}) \to 0.$$ This follows
from \fullref{prop:local.coefficients}.

The next assertion is that certain relations are satisfied in the
Lie algebra $gr^*(P\varSigma_n)$: If $\{i,j\} \cap \{s,t\} = \emptyset$,
then $[\chi_{j,i}, \chi_{s,t}] = 0$ by one of McCool's relations in
\fullref{thm:McCool relations}.

Next notice that one of the horizontal $4T$ relations follows
directly from McCool's identity in \fullref{thm:McCool
relations}. Since $[ \chi_{i,j} \cdot \chi_{k,j}, \chi_{i,k}] = 1$
for $i,j,k$ distinct, it follows that $$[\chi_{i,k}, \chi_{i,j} +
\chi_{k,j}] = 0$$ on the level of Lie algebras.

In addition, $[ \chi_{k,j}, \chi_{s,j}] = 0$ on the level of Lie
algebras if $\{s,k\} \cap \{j \}= \emptyset$ by inspection of \fullref{thm:McCool relations}.

That $[\chi_{k,i}, \chi_{j,i} + \chi_{j,k}]$ is non-zero follows
from \cite{cp}. The details are omitted: they are a direct
computation using the Johnson homomorphism together with structure
for the Lie algebra of derivations of a free Lie algebra.

Next notice that $H_1(P\varSigma_{n}) = \oplus_{\binom{n}{2}}\mathbb
Z$. Project it to the summand with basis $\chi_{i,n}$ for $ 1 \leq i
\leq n-1$. Denote by $\gamma$ the following  composition 
$$K_n\to P\varSigma_{n}\to H_1(P\varSigma_{n}) \to
\oplus_{n-1} \mathbb Z.$$ 
It is evidently a split epimorphism as
$[\chi_{i,n},\chi_{j,n}] = 1.$ The remaining properties follow by
inspection. \fullref{thm:PSigma.n.Lie.algebras} follows.

\section[{Proof of \ref{thm:PSigma_{n,+}}}]{Proof of \fullref{thm:PSigma_{n,+}}}

\fullref{prop:local.coefficients} gives that the action of
$P\varSigma_{n-1}^+$ on $H_1(K_n^+)$ is trivial. There are two
consequences of this fact.

The first consequence is that there is a split short exact sequence
of Lie algebras 
$$0 \to\ gr^*(K_n^+) \to\ gr^*(P\varSigma_{n}^+) \to\
gr^*(P\varSigma_{n-1}^+) \to\ 0$$ 
by \fullref{prop:local.coefficients} or \cite{k,fr}.

Notice that one of the horizontal $4T$ relations follows directly
from McCool's identity in \fullref{thm:McCool relations}. As
checked in the proof of \fullref{thm:PSigma.n.Lie.algebras}, the
relation $$[\chi_{i,k},\chi_{i,j}+ \chi_{k,j}] = 0$$ is satisfied on
the level of Lie algebras.

The remaining two relations $[\chi_{j,k}, \chi_{s,t}] = 0$ if
$\{j,k\} \cap \{s,t \}= \emptyset$, and $[\chi_{j,k}, \chi_{s,j}] = 0$ if
$\{s,k\} \cap \{j\}= \emptyset$ follow by inspection of McCool's
relations. Thus the asserted structure of Lie algebra follows.

The second consequence is that the local coefficient system for the
Lyndon--Hochschild--Serre spectral sequence of the extension
$$1 \to K_n^+ \to P\varSigma_{n}^+ \to P\varSigma_{n-1}^+ \to 1$$
has trivial local coefficients in cohomology. Since $K_n^+$ is a
free group, it has torsion free cohomology which is concentrated in
degrees at most $1$. Thus the $E_2$--term of the spectral sequence
splits as a tensor product
$$H^*(P\varSigma_{n-1}^+) \otimes_{\mathbb Z} H^*( K_n^+).$$ Since the extension is
split, all differentials are zero, and the spectral sequence
collapses at the $E_2$--term.

Note that $P\varSigma_2^+$ is isomorphic to the integers. Thus the
cohomology of $P\varSigma_n^+$ is torsion free and a basis for the
cohomology is given as stated in the theorem by inspection of the
$E_2$--term of the spectral sequence by induction on $n$.

It thus remains to work out the product structure in cohomology
which is asserted to be
\begin{itemize}
     \item ${\chi_{k,i}^*}^2 = 0$ for all $i > k$, and
      \item $\chi_{i,j}^*[\chi_{i,k}^*- \chi_{j,k}^*]=0$ for $k < j < i$,
\end{itemize} 
for which $\chi_{i,k}^*$ is the basis
element in cohomology dual to $\chi_{i,k}$.

Notice that ${\chi_{i,k}^*}^2 = 0$. It suffices to work out the
relation $\chi_{3,2}^*\cdot [\chi_{3,1}^*- \chi_{2,1}^*]=0$ in case
$n = 3$ by the natural projection maps.

Consider the natural split epimorphism obtained by restriction to
$$\pi|_{P\varSigma_3^+}\co P\varSigma_3^+\to P\varSigma_2^+.$$ The kernel 
is a free group on two
letters $\chi_{3,1}, \chi_{3,2}$ and the spectral sequence of the
extension collapses. On the level of cohomology, it suffices to work
out the value of $\chi_{3,1}^*\cdot \chi_{3,2}^*$. Thus the product
$\chi_{3,1}^*\cdot \chi_{3,2}^*$ is equal to the linear combination
$A\chi_{2,1}^*\cdot \chi_{3,1}^* + B\chi_{2,1}^*\cdot \chi_{3,2}^*$
for scalars $A$ and $B$.

Next, consider McCool's relations (as stated in \ref{thm:McCool
relations}) which gives
$$\chi_{3,1}\cdot \chi_{2,1}\cdot \chi_{3,2} = 
\chi_{3,2}\cdot \chi_{3,1}\cdot \chi_{2,1}.$$ Thus
the commutator $[\chi_{3,1}\cdot \chi_{2,1},\chi_{3,2}]$ is $1$ in
$P\varSigma_3$. Since this commutator is trivial, there is an
induced homomorphism
$$\rho\co  \mathbb Z \times \mathbb Z \to P\varSigma_3^+$$ defined by the
equation
     \[
    \rho((a,b))=
    \begin{cases}
     \chi_{3,2} & \text{if $(a,b) = (1,0)$ and} \\
    \chi_{3,1}\cdot \chi_{2,1} & \text{if $(a,b) = (0,1)$.}
    \end{cases}
    \]
Let $(1,0)^*$ and $(0,1)^*$ denote the two associated natural
classes in $H^1(\mathbb Z \times \mathbb Z; \mathbb Z)$ which are
dual to the natural homology basis. Notice that

\begin{itemize}
    \item $\rho^*(\chi_{2,1}^*) = (0,1)^*$,
    \item $\rho^*(\chi_{3,2}^*) = (1,0)^*$ and
    \item $\rho^*(\chi_{3,1}^*) = (0,1)^*$.
\end{itemize}

Let $(1,1)^*$ denote the fundamental cycle of $H^2(\mathbb Z \times
\mathbb Z)$ given by the cup product $(1,0)^*\cdot(0,1)^*$. Consider
$$-(1,1)^* = \rho^*(\chi_{3,1}^*\cdot \chi_{3,2}^*) = \rho^*(A\chi_{2,1}^*\cdot
\chi_{3,1}^* + B\chi_{2,1}^*\cdot \chi_{3,2}^*) = -B(1,1)^*.$$ Thus
$B = 1$.

Next, consider the automorphisms of $P\varSigma_3^+$. Recall that
there are automorphisms of $P\varSigma_n^+$ specified by the
conjugation action of $\tau_i$ on the elements $\chi_{s,t}$ given by
the formulas:
    \[
    \tau_i(\chi_{s,t})\tau_i^{-1}=
    \begin{cases}
    \chi_{s,t} & \text{if $\{i\} \cap \{s,t\} = \emptyset$},\\
    \chi_{s,i}^{-1} & \text{ if $t = i$},\\
    \chi_{i,t} & \text{if $s = i$.}
    \end{cases}
    \]Notice that $\tau_1$ leaves $P\varSigma_3^+$ invariant and is thus an
    automorphism of $P\varSigma_3^+$. Thus on the level of $H_1(P\varSigma_3^+)$,
    conjugation by $\tau_1$ denoted $\phi_1$, is given by the formula
\begin{enumerate}
    \item ${\phi_1}_*(\chi_{2,1}) = -\chi_{2,1}$,
    \item ${\phi_1}_*(\chi_{3,1}) = -\chi_{3,1}$ and
    \item ${\phi_1}_*(\chi_{3,2}) = \chi_{3,2}$.
\end{enumerate}

Thus on the level of cohomology,
\begin{enumerate}
    \item $\phi_1^*(\chi_{2,1}^*) = -\chi_{2,1}^*$,
    \item $\phi_1^*(\chi_{3,1}^*) = -\chi_{3,1}^*$ and
    \item $\phi_1^*(\chi_{3,2}^*) = \chi_{3,2}^*$.
\end{enumerate}

Apply the automorphism $\phi_1^*$ to the equation $\chi_{3,1}^*\cdot
\chi_{3,2}^*= A\chi_{2,1}^*\cdot \chi_{3,1}^* + B\chi_{2,1}^*\cdot
\chi_{3,2}^*$ where $B=1$ to obtain the following.

\begin{itemize}
    \item $-\chi_{3,1}^*\cdot \chi_{3,2}^* = A \chi_{2,1}^*\cdot \chi_{3,1}^* -
B\chi_{2,1}^*\cdot \chi_{3,2}^*$.
    \item $A \chi_{2,1}^*\cdot \chi_{3,1}^* - B\chi_{2,1}^*\cdot \chi_{3,2}^* = -[A
\chi_{2,1}^*\cdot \chi_{3,1}^* + B\chi_{2,1}^*\cdot \chi_{3,2}^*]$
in a free abelian group of rank two with basis $\{\chi_{2,1}^*\cdot
\chi_{3,1}^*, \chi_{2,1}^*\cdot \chi_{3,2}^*\}$.
    \item Hence $A = 0$ and $\chi_{3,2}^*\cdot [\chi_{3,1}^*- \chi_{2,1}^*]=0$.
    \item It follows that $\chi_{i,j}^*\cdot [\chi_{i,k}^*- \chi_{j,k}^*]=0$ 
for $k < j < i $
    by a similar argument.

\end{itemize} The Theorem follows.

\section[{Proof of \ref{thm:A(n)}}]{Proof of \fullref{thm:A(n)}}

Recall that $H_1(IA_n)$ is a free abelian group of rank
$n\binom{n}{2}$ with basis given by the equivalence classes of
$\chi_{i,k}$ for $ i \neq k$ with $1 \leq i,k \leq n$, and
$\theta(k;[s,t])$ for $k,s,t$ distinct positive natural numbers with
$s< t$ \cite{ka,f,cp,sa}. Thus the natural composite
\[
\begin{CD}
P\varSigma_{n}^+ @>{}>> IA_n @>{}>> IA_n/[IA_n, IA_n]
\end{CD}
\] is a split monomorphism on the level of the first homology
group with image spanned by $\chi_{k,j}$ for $k > j$.

The cohomology algebra of $P\varSigma_n^+$ is generated as an
algebra by elements of degree $1$ given by $\chi_{k,j,}^*$ for $k >
j$ in the dual basis for $H_1(P\Sigma_{n}^+)$ by \fullref{thm:PSigma_{n,+}}. Thus the composite $P\Sigma_{n}^+ \to IA_n
\to IA_n/[IA_n, IA_n]$ is a split surjection in integral cohomology.
That the image is not $GL(n, \mathbb Z)$--invariant follows by an
inspection. Since the composite
\[
\begin{CD}
H^*(IA_n/[IA_n,IA_n]) @>{A^*}>> H^*(IA_n) @>{}>> H^*(P\Sigma_{n}^+)
\end{CD}
\] is an epimorphism, and $A^*$ is $GL(n,\mathbb Z)$--equivariant,
the image of $A^*$ contains $\mathbb M_n$.

Notice that the classifying space $BIA_n/[IA_n, IA_n]$ is homotopy
equivalent to a product of circles $(S^1)^{n{\binom{n}{2}}}.$
Furthermore, the composite $P\varSigma_{n}^+ \to IA_n \to
IA_n/[IA_n, IA_n]$ induces a split epimorphism on integral
cohomology. That the map $$BP\varSigma_{n}^+ \to
(S^1)^{n{\binom{n}{2}}}$$ is split after suspending once follows
directly from the next standard property of maps into products of
spheres.

\begin{lem} \label{lem:maps.to.products.of.spheres}
Let $f\co X \to Y$ be a continuous map which satisfies the following
properties.
\begin{enumerate}
  \item The space $X$ is homotopy equivalent to a CW--complex,
  \item $Y$ is a finite product of spheres of dimension at least $1$
  and
  \item the map $f$ induces a split epimorphism on integral
  cohomology.
\end{enumerate}
Then
\begin{enumerate}
  \item the suspension of $X$, $\Sigma(X)$, is a retract of  $\Sigma(Y)$
  \item $\Sigma(X)$ is homotopy equivalent to a bouquet of
  spheres and
  \item the Freudenthal suspension $E\co X \to \Omega \Sigma(X)$ factors through a
map $Y \to \Omega \Sigma(X)$.
\end{enumerate}
\end{lem}

Hence there is a map $$\Theta\co  BIA_n \to \Omega
\Sigma(P\Sigma_{n}^+)$$ which gives a factorization of the
Freudenthal suspension. The Theorem follows.

\section[{Proof of \ref{thm:braid.permutation}}]{Proof of \fullref{thm:braid.permutation}}

The purpose of this section is to give the proof of \fullref{thm:braid.permutation} along with other information. There is a
homomorphism $$\rho_n\co \mathrm{BP}_n \to \Sigma_n$$ defined by
$\rho(\sigma_i)=\xi_i = \rho(\xi_i)$. One way to see that such a
homomorphism exists is to consider the pullback diagram
\[
\begin{CD}
\Xi_n @>{\pi}>> \Sigma_n    \\
  @V{i}VV          @VV{}V \\
\Aut_n @>{}>> GL(n,\mathbb Z)
\end{CD}
\] and to observe that $\mathrm{BP}_n$ is a subgroup of $\Xi_n$.

Recall that the group $\mathrm{BP}_n$ is characterized as the
subgroup of $\Aut(F_n)$ consisting of automorphisms $\phi\in
\Aut(F_n)$ of {\it permutation--conjugacy type} which satisfy
\begin{equation}
%%\label{b_p_type}
(x_i)\phi=w_i^{-1}x_{\lambda(i)}w_i
\end{equation}
for some word $w_i\in F_n$ and permutation $\lambda \in \Sigma_n$
the symmetric group on $n$ letters (where the action of an element
$\phi$ in $\Aut(F_n)$ is from the right) \cite{FRR2}. Thus

\begin{equation}\label{b_p_type}
((x_i)\phi \lambda^{-1})=
((w_i^{-1})\lambda^{-1})(x_{i})((w_i)\lambda^{-1}).
\end{equation} Furthermore $\sigma_i\circ \xi_i =
\chi_{i+1,i}$ and so $\sigma_i = \chi_{i+1,i} \circ \xi_i^{-1}.$ The
group $\mathrm{BP}_n$ is thus generated by (1) $\Sigma_n$ and (2)
$P\varSigma_n$. Thus the kernel of the natural map $\rho_n\co
\mathrm{BP}_n \to \Sigma_n$ is $P\varSigma_n$. \fullref{thm:braid.permutation} follows.

\section{Problems}
\begin{enumerate}
\item Is $P\Sigma_{n}^+$ isomorphic to $P_n$ ? If $n=2,3$, the
answer is clearly yes. Note added after this paper was submitted:
Bardakov \cite{BM} has shown that if $n > 3$, the groups $P\Sigma_{n}^+$ and
$P_n$ are not isomorphic.

\item Let $\mathcal K_n$ denote the subgroup of $IA_n$ generated
by all of the elements $\theta(k;[s,t])$. Identify the structure of
$gr^*(\mathcal K_n)$.
\item Are the natural maps $$gr^*(P\Sigma_{n}^+)\to gr^*(IA_n),$$
$$gr^*(P\Sigma_{n})\to gr^*(IA_n)$$ and/or
$$gr^*(\mathcal K_n)\to gr^*(IA_n)$$ monomorphisms ?
\item Give the Euler--Poincar\'e series for $U[gr^*(IA_n)\otimes \mathbb Q]$ and
$U[gr^*(P\Sigma_n)\otimes \mathbb Q]$ where $U[\mathcal L]$ denotes
the universal enveloping algebra of a Lie algebra $\mathcal L$.
\item Give the Euler--Poincar\'e series for $\mathbb M_n$.
\item Observe that the
kernel of the natural map out of the free product
$P\Sigma_{n}*\mathcal K_n \to IA_n$ is a surjection with kernel a
free group.  The natural morphism of Lie algebras
$$\iota\co gr^*(P\Sigma_{n}) \amalg gr^*(\mathcal K_n)\to gr^*(IA_n)$$
is an epimorphism (where $\amalg$ denotes the free product of Lie
algebras). Is the kernel a free Lie algebra ?

\end{enumerate}

\section{Acknowledgments}
This work was started during the visit of the first and the third
authors to the National University of Singapore. They are very
thankful for hospitality to the Mathematical Department of this
university and especially to Jon Berrick.

The first author was partially supported by National Science
Foundation under Grant No. 0305094 and the Institute for Advanced
Study. The third author was supported in part by INTAS grant
No~03-5-3251 and the ACI project ACI-NIM-2004-243 ``Braids and
Knots". The fourth author was partially supported by a grant from
the National University of Singapore. The first and third authors
were supported in part by a joint CNRS-NSF grant No~17149.

\bibliographystyle{gtart}
\bibliography{link}

\begin{thebibliography}{}
\providecommand\bibmarginpar{\leavevmode\marginpar}
\def\urlstyle#1{{\tt #1}}

\bibitem{bardakov}
\textbf{V\,G Bardakov}, \href{http://dx.doi.org/10.1023/A:1025913505208}
  {\emph{The structure of a group of conjugating automorphisms}}, Algebra
  Logika 42 (2003) 515--541, 636 \xox{MR}{2025714}

\bibitem{BM}
\textbf{V\,G Bardakov}, \textbf{R Mikhailov}, \emph{On certain questions of the
  free group automorphisms theory} \xox{arXiv}{math.GR/0701441}

\bibitem{ccp}
\textbf{D\,C Cohen}, \textbf{F\,R Cohen}, \textbf{S Prassidis},
  \emph{Centralizers of {L}ie algebras associated to descending central series
  of certain poly-free groups}, J. Lie Theory 17 (2007) 379--397
  \xox{MR}{2325705}

\bibitem{cp}
\textbf{F\,R Cohen}, \textbf{J Pakianathan}, \emph{Notes on automorphism
  groups}, in preparation

\bibitem{fr}
\textbf{M Falk}, \textbf{R Randell},
  \href{http://dx.doi.org/10.1007/BF01394779} {\emph{The lower central series
  of a fiber-type arrangement}}, Invent. Math. 82 (1985) 77--88
  \xox{MR}{808110}

\bibitem{f}
\textbf{B Farb}, \emph{The Johnson homomorphism for {\rm Aut(}$F_n${\rm )}}, in
  preparation

\bibitem{FRR2}
\textbf{R Fenn}, \textbf{R Rim{\'a}nyi}, \textbf{C Rourke},
  \href{http://dx.doi.org/10.1016/0040-9383(95)00072-0} {\emph{The
  braid-permutation group}}, Topology 36 (1997) 123--135 \xox{MR}{1410467}

\bibitem{g}
\textbf{G Galatius}, \emph{Stable homology of automorphism groups of free
  groups} \xox{arXiv}{math.AT/0610216}

\bibitem{dg}
\textbf{D\,L Goldsmith},
  \href{http://projecteuclid.org/getRecord?id=euclid.mmj/1029002454} {\emph{The
  theory of motion groups}}, Michigan Math. J. 28 (1981) 3--17 \xox{MR}{600411}

\bibitem{h}
\textbf{R Hain}, \href{http://dx.doi.org/10.1090/S0894-0347-97-00235-X}
  {\emph{Infinitesimal presentations of the {T}orelli groups}}, J. Amer. Math.
  Soc. 10 (1997) 597--651 \xox{MR}{1431828}

\bibitem{y}
\textbf{Y Ihara}, \emph{Braids, {G}alois groups, and some arithmetic
  functions}, from: ``Proceedings of the International Congress of
  Mathematicians, Vol.\ I, II (Kyoto, 1990)'', Math. Soc. Japan, Tokyo (1991)
  99--120 \xox{MR}{1159208}

\bibitem{jmm}
\textbf{C Jensen}, \textbf{J McCammond}, \textbf{J Meier},
  \href{http://dx.doi.org/10.2140/gt.2006.10.759} {\emph{The integral
  cohomology of the group of loops}}, Geom. Topol. 10 (2006) 759--784
  \xox{MR}{2240905}

\bibitem{j}
\textbf{D Johnson}, \href{http://dx.doi.org/10.1007/BF01363897} {\emph{An
  abelian quotient of the mapping class group $\mathcal{I}_g$}}, Math. Ann. 249
  (1980) 225--242 \xox{MR}{579103}

\bibitem{ka}
\textbf{N Kawazumi}, \emph{Cohomological aspects of Magnus expansions}
  \xox{arXiv}{math.GT/0505497}

\bibitem{k}
\textbf{T Kohno}, \href{http://dx.doi.org/10.1007/BF01394779} {\emph{S\'erie de
  {P}oincar\'e-{K}oszul associ\'ee aux groupes de tresses pures}}, Invent.
  Math. 82 (1985) 57--75 \xox{MR}{808109}

\bibitem{magnus}
\textbf{W Magnus}, \href{http://dx.doi.org/10.1007/BF02545673} {\emph{\"{U}ber
  {$n$}-dimensionale {G}ittertransformationen}}, Acta Math. 64 (1935) 353--367
  \xox{MR}{1555401}

\bibitem{mks}
\textbf{W Magnus}, \textbf{A Karrass}, \textbf{D Solitar}, \emph{Combinatorial
  group theory}, second edition, Dover Publications Inc., Mineola, NY (2004)
  \xox{MR}{2109550}Presentations of groups in terms of generators and relations

\bibitem{mc}
\textbf{J McCool}, \emph{On basis-conjugating automorphisms of free groups},
  Canad. J. Math. 38 (1986) 1525--1529 \xox{MR}{873421}

\bibitem{m}
\textbf{S Morita}, \href{http://dx.doi.org/10.1215/S0012-7094-93-07017-2}
  {\emph{Abelian quotients of subgroups of the mapping class group of
  surfaces}}, Duke Math. J. 70 (1993) 699--726 \xox{MR}{1224104}

\bibitem{N}
\textbf{J Nielsen}, \href{http://dx.doi.org/10.1007/BF01458209} {\emph{\"{U}ber
  die {I}somorphismen unendlicher {G}ruppen ohne {R}elation}}, Math. Ann. 79
  (1918) 269--272 \xox{MR}{1511927}

\bibitem{pe}
\textbf{A Pettet}, \href{http://dx.doi.org/10.2140/agt.2005.5.725} {\emph{The
  {J}ohnson homomorphism and the second cohomology of {${\rm IA}\sb n$}}},
  Algebr. Geom. Topol. 5 (2005) 725--740 \xox{MR}{2153110}

\bibitem{s}
\textbf{T Sakasai}, \href{http://dx.doi.org/10.1016/j.topol.2004.08.002}
  {\emph{The {J}ohnson homomorphism and the third rational cohomology group of
  the {T}orelli group}}, Topology Appl. 148 (2005) 83--111 \xox{MR}{2118957}

\bibitem{sa}
\textbf{T Satoh}, \href{http://dx.doi.org/10.1017/S0305004106009959} {\emph{The
  abelianization of the congruence {IA}-automorphism group of a free group}},
  Math. Proc. Cambridge Philos. Soc. 142 (2007) 239--248 \xox{MR}{2314598}

\end{thebibliography}

\end{document}